\documentclass{cimsart}                    %% Input Standard CPAM class.
\received{Month 199X}       %       \revised{date of revision}
\yearofpublication{199X}                \volume{000}
\startingpage{1}                        \cccline{00}

\authorheadline{O. Costin and M. D. Kruskal}
\titleheadline{Integrability of difference
  equations}
\ifx\pdftexversion\undefined
  \usepackage[dvips]{graphics}
\else
  \usepackage[pdftex]{graphics}
\fi
\usepackage{amsmath} \usepackage{latexsym} \usepackage{amssymb}
\newtheorem{Lemma}[theorem]{Lemma}

\newtheorem{Proposition}[theorem]{Proposition}
\newtheorem{Definition}[theorem]{Definition}

\newtheorem{Note}[theorem]{Note}
\newtheorem{Remark}[theorem]{Remark}
\begin{document}                        %% Standard LaTeX command

%%      ---------------------------------------------------------------------
%%      -------------------------------- TITLE ------------------------------
%%      ---------------------------------------------------------------------

%%      Please replace material within asterisks, including the asterisks
%%      themselves, with the correct information.

\title{Analytic methods for obstruction to integrability in discrete
  dynamical systems} \gdef\shorttitle{Integrability of difference
  equations}              %% Enter the title here.

%%      ---------------------------------------------------------------------
%%      ------------------------------ AUTHORS ------------------------------
%%      ---------------------------------------------------------------------

%%      Please replace material within asterisks with correct information.

\author{O  Costin \\            %% Enter (first) author's name,
\affil  Rutgers University
\\ Piscataway, NJ 08854-8019\\  
% \\ AND \\
% \\ M  D  Kruskal
% \affil Rutgers University
\\ Piscataway, NJ 08854-8019}

%%      ---------------------------------------------------------------------
%%      --------------------------- DEDICATION  (OPTIONAL)------------------- 
%%      ---------------------------------------------------------------------

%       Uncomment the following line to insert a dedication.

%\dedication{* Dedication *}            %% Enter dedication between braces.

%%      ---------------------------------------------------------------------
\maketitle   
%%      ---------------------------------------------------------------------
%%      --------------------------- ABSTRACT (OPTIONAL)----------------------
%%      ---------------------------------------------------------------------

%% ***** UNCOMMENT THE FOLLOWING AND INSERT YOUR ABSTRACT IF YOU HAVE ONE
\def\z{\noindent}
\def\CC{\mathbb{C}}
\def\RR{\mathbb{R}}
\def\NN{\mathbb{N}}
\def\ZZ{\mathbb{Z}}
\def\lap{\mathcal{L}}
\def\bor{\mathcal{B}}
\def\phi{\varphi}
\def\bfk{\mathbf{k}}
 \begin{abstract}
 
  A unique analytic continuation result is proved for solutions of a
  relatively general class of difference equations, using techniques
  of generalized Borel summability.
  
  This continuation allows for Painlev\'e property methods to be
  extended to difference equations.
 
  It is shown that the Painlev\'e property (PP) induces, under
  relatively general assumptions, a dichotomy within first order
  difference equations: all equations with PP can be solved in closed
  form; on the contrary, absence of PP implies, under some further
  assumptions, that the local conserved quantities are strictly local
  in the sense that they develop singularity barriers on the boundary
  of some compact set.
 
  The technique produces analytic formulas to describe fractal sets
  originating in polynomial iterations.

 \end{abstract}

%%      ---------------------------------------------------------------------
%%      ------------------- TABLE OF CONTENTS (OPTIONAL) --------------------
%%      ---------------------------------------------------------------------

%% ***** IF YOUR PAPER IS OVER 40 PAGES AND YOU WISH TO HAVE A TABLE
%% ***** OF CONTENTS, PLEASE UNCOMMENT THE FOLLOWING LINE

\tableofcontents
\section{Introduction and main results}
Solvability of difference equations as well as chaotic behavior have
stimulated extensive research. For differential equations the
Painlev\'e test, which consists in checking whether all solutions of a
given equation are free of movable non-isolated singularities provides
a convenient and effective tool in detecting integrable cases (see
\S\ref{hist1}).

A difficulty in applying Painlev\'e's methods to difference equations
resides in extending the solutions, which are defined on a discrete
set, to the complex plane of the independent variable in a natural and
effective fashion, when, in the interesting cases, there is no
explicit formula for them.  A number of alternative approaches, but no
genuine analog of the Painlev\'e test, have been proposed, see
\cite{[AHH]} \cite{[CM]} \cite{Joshi} \cite{RG} (a comparative
discussion of the various approaches is presented in \cite{[AHH]}).

The present paper proposes a natural way, based on generalized Borel
summabillity, to extend the solutions in the complex plane
(Theorem~\ref{T0} below), allowing for a definition of a discrete
Painlev\'e test.  Subsequent analysis shows that the test is sharp in
a class of first order difference equations: those passing the test are
explicitly solvable (Theorem~\ref{T1}) while polynomial equations
failing the test exhibit chaotic behavior and their local conserved
quantities (see \S\ref{S19}) develop barriers of singularities along
fractal sets (Theorem~\ref{Barrier}).

The approach also allows for a detailed study of analytic properties
near these singularity barriers as well as finding rapidly convergent
series representing the corresponding fractal curves
(Theorem~\ref{T5}).
\subsection{Setting}\label{1.1}
\z We consider difference systems of equations which can be brought to
the form
\begin{equation}
  \label{eq:Bra1}
  \mathbf{x}(n+1)=\hat\Lambda\left(I+ \frac{1}{n}\hat{A}\right)\mathbf{x}(n) +\mathbf{g}(n,\mathbf{x}(n))
\end{equation}
where $\hat\Lambda $ and $\hat{A} $ are constant coefficient matrices, $\mathbf{g}$ is
convergently given for small $\mathbf{x}$ by
\begin{equation}
  \label{eq:defg}
\mathbf{g}(n,\mathbf{x})=\sum_{\mathbf{k}\in\NN^m}
\mathbf{g}_\mathbf{k}(n)\mathbf{x}^\mathbf{k} 
\end{equation}
with $\mathbf{g}_\mathbf{k}(n)$ analytic in $n$ at infinity and
\begin{equation}
  \label{eq:condg}
  \mathbf{g}_\mathbf{k}(n)=O(n^{-2}) \ \mbox{as }n\rightarrow\infty,\
  \mbox{if }\sum_{j=1}^m k_j\le 1
\end{equation}
 under nonresonance conditions: Let
$\boldsymbol{\mu}=(\mu_1,...,\mu_n)$ and $\mathbf a=(a_1,...,a_n)$ where
$e^{-\mu_k}$ are the eigenvalues of $\hat{\Lambda}$ and the $a_k$ are
the eigenvalues of $\hat{A}$.  Then the nonresonance condition is
\begin{equation}
  \label{eq:nonres}
  \left(\mathbf{k}\cdot \boldsymbol\mu=0 \mod 2\pi i \ \ \mbox{with}\ \ \mathbf{k} \in\ZZ^{m_1}
  \right)
\Leftrightarrow \mathbf{k}=0.
\end{equation}
We consider the solutions of (\ref{eq:Bra1}) which are small as $n$
becomes large.
\subsection{Analyzability: transseries and generalized Borel summability} These concepts
were introduced by \'Ecalle in the fundamental work \cite{Eca84}.
Analyzability of difference equations was shown in
\cite{Braaksma,Eca84}.  We give below a brief description of the
concepts effectively used in the present paper and refer to
\cite{Duke,Braaksma} for a general theory.  An
expression of the form
\begin{equation}
  \label{eq:transszf}
  \tilde{{\bf x}}(t)
:=\sum_{{\bf k}\in\NN^{m}}{\bf C}^{\bf k} e^{-{\bf k}\cdot
  \boldsymbol{\mu} t}
t^{{\bf k}\cdot {\bf a}} \tilde{{\bf x}}_\bfk(t)
\end{equation}
where $ \tilde{{\bf x}}_\bfk(t)$ are formal power series in powers of
$t^{-1}$ is an exponential power series; it is a transseries as $t\to
+\infty$ if $\Re(\mu_j)>0$ for all $j$ with $1\le j\le m$.  Such a transseries
is Borel summable as $t\to
+\infty$ if there exist constants $A,\nu>0$ and a family of functions
  \begin{eqnarray}
    \label{estX}&
\text{ ${\bf
  X}_{\mathbf{k}}$ analytic in a sectorial neighborhood $\mathcal{S}$
  of $\RR^+$, satisfying}\nonumber\\&
    \displaystyle \sup_{p\in\mathcal{S},\bfk\in\NN^{m}}\left|A^{|\bf k|} e^{-\nu|p|}
\mathbf{X}_{\mathbf{k}}\right|<\infty
  \end{eqnarray}
such that the functions ${\bf x}_{\mathbf{k}}$ defined by
\begin{equation}
  \label{eq:Bory} {\bf x}_{\mathbf{k}}(t)=\int_0^{\infty}e^{-tp}{\bf
  X}_{\mathbf{k}}(p)dp
\end{equation}
are asymptotic to the series $\tilde{\mathbf{x}}_{\mathbf{k}}$ i.e.
\begin{equation}
  \label{asyy1}
  {\bf x}_{\mathbf{k}}(t)\ \sim\ \tilde{\bf x}_{\bf k}(t)\ \ (t\to+\infty)
\end{equation}
It is then easy to check that condition (\ref{estX}) implies that the
sum
\begin{equation}
  \label{eq:transsS}
  {\bf x}(t)
=\sum_{{\bf k}\in\NN^{n_0}}{\bf C}^{\bf k} e^{-{\bf k}\cdot
  \boldsymbol{\mu} t}
t^{{\bf k}\cdot {\bf a}}{\bf x}_\bfk(t)
\end{equation}
is convergent in the half plane $\mathbb{H}=\{t:\Re(t)>t_0\}$, for
$t_0$ large enough.  The function ${\bf x}$ in (\ref{eq:transsS}) is
by definition the Borel sum of the transseries $\tilde{\bf x}$ in
(\ref{eq:transszf}). Generalized Borel summability allows for
singularities of $\mathbf{X}_{\bf k}$ of certain types along $\RR^+$.
The transseries $\tilde{\bf x}$ is (generalized) Borel summable in the
direction $e^{i\phi}\RR^+$ if $\tilde{\bf x}(\cdot e^{-i\phi})$ is
(generalized) Borel summable. (Generalized) Borel summation is known
to be an extended isomorphism between transseries and their sums, see
\cite{Eca84}, \cite{EcalleNATO}, \cite{Duke}.
\subsubsection{Transseries for difference equations}\label{resB} Braaksma \cite{Braaksma} showed that the recurrences (\ref{eq:Bra1}) 
 posess $l$-parameter transseries solutions of the form
 (\ref{eq:transszf}) with $t=n$ where $ \tilde{{\bf x}}_\bfk(n)$ are
 formal power series in powers of $n^{-1}$ and $l\le m$ is chosen such
 that, after reordering the indices, we have $\Re(\mu_j)>0$ for $1\le
 j\le l$.
 
 It is shown in \cite{Braaksma} and \cite{Kuik} that these transseries
 are generalized Borel summable in any direction and Borel summable in
 all except $m$ of them and that
\begin{equation}
  \label{eq:transsSn}
  {\bf x}(n)
=\sum_{{\bf k}\in\NN^{l}}{\bf C}^{\bf k} e^{-{\bf k}\cdot
  \boldsymbol{\mu} n}
n^{{\bf k}\cdot {\bf a}}{\bf x}_\bfk(n)
\end{equation} 
is a solution of (\ref{eq:Bra1}), if $n>y_0$, $t_0$ large enough.
 \subsection{Uniqueness of  continuation from $\NN$ to $\CC$}
The values of $\bf x$ on the integers uniquely determine $\bf x$.
\begin{theorem}\label{T0}  In the assumptions in \S\ref{1.1} and  \ref{resB},
  define the continuation of ${\bf x}_\bfk(n)$ in the half plane
  $\{t:\Re(t)>t_0\}$ by ${\bf x}(t)$, cf.  (\ref{estX})--(\ref{eq:transsS}).
  
  The following uniqueness property holds. If in the assumptions
  (\ref{estX})--(\ref{eq:transsS}) we have ${\bf x}(n)=0$ for all
  except possibly finitely many $n\in\NN$, then $ {\bf x}(t)$ =0 for
  all $t\in\CC,\,\Re(t)>t_0$.
\end{theorem}
\z The proof is given in \S\ref{Uniq}.
\subsection{Continuation of solutions of difference equations to the complex $n$ plane}
The representation (\ref{eq:transsSn}) and Theorem~\ref{T0} make the following
definition natural.
\subsection{Continuability and singularities} 
The function ${\bf x}$ is analytic in $\mathbb{H}$ and has, in
general, nontrivial singularities in $\CC\setminus\mathbb{H}$. The
results in \cite{Invent}, extended to difference equations in
\cite{Braaksma, Braaksma-Kuik,Kuik}, give constructive methods to
determine those singularities that arise near the boundary of
$\mathbb{H}$; these form, generically, nearly periodic arrays.
\subsection{Integrability}\label{Int} In particular,
Painlev\'e's test of integrability (absence of movable non-isolated
singularites) extends then to difference equations.

As in the case of differential equations, fixed singularities are singular
points whose location is the same for all solutions; they define a common
Riemann surface. Other singularities (i.e., whose location depends on initial
data) are called {\em movable}.
 
\begin{Definition}\label{DPP}
  We say that a difference equation has the Painlev\'e
  property if its solutions are analyzable and their analytic continuations on
  a Riemann surface {\em common to all solutions}, have only {\em isolated}
  singularities.
\end{Definition}

\z {\bf Note.} We follow the usual convention that an isolated
singular point of an analytic function $f$ is a point $z_0$ such that
$f$ is analytic in some disk centered at $z_0$ except perhaps at $z_0$
itself.  Branch points are thus not isolated singularities and neither
are singularity barriers; it is worth noting, however, that for
differential equations there exist equations sometimes considered
integrable (the {\em Chazy equation}, a third order nonlinear one is
the simplest known example) whose solutions exhibit singularity
barriers.
\subsection{First order autonomous equations}
\label{1.5}
These are equations of the type
\begin{equation}
  \label{eq:01}
  x_{n+1}=G(x_n):=ax_n+F(x_n)
\end{equation}
Some analyticity assumptions on $F$ are required for our method to
apply. We define a class of single valued functions closed under all
algebraic operations and composition (the latter is needed since $x_n$
written in terms of $x_0$ involves repeated composition).

We need to allow for singular behavior in $F$, and meromorphic
functions are obviously not closed under composition. The following
definition formalizes an extension of meromorphic functions, often
used informally in the theory of integrability.
\begin{Definition}
  We define the "mostly analytic functions" to be the class
    $\mathcal{M}$ of functions analytic in the complement of a closed
    countable set (which may depend on the function).
\end{Definition}

\begin{Lemma}\label{L22}
  (a) The class $\mathcal{M}$ is closed under addition, multiplication
  and multiplication by scalars, and also under division and
  composition between (nonconstant) functions. It includes meromorphic
  functions.

 (b) If $G\in \mathcal{M}$ is not a constant, then the equation $G(x)=y$
  has solutions for all large enough $y$.

(c). The class $\mathcal{M}_0$ of $G\in \mathcal{M}$, with $G$ analytic at zero, $G(0)=0$
 and $0<|G'(0)|<1$ is closed under composition.

In particular, $G^{\circ m}\in \mathcal{M}$ for $m\ge 1$.
\end{Lemma}
\z {\em Proof}. All properties in (a) are obvious except for closure
under composition and division, proved in \S~\ref{PL4}; (b) follows
from the proof of Lemma~\ref{L12}.  (c) is easily shown using (a).
\subsection{Classification of equations of type (\ref{eq:01}) with respect to integrability}
\begin{theorem}\label{T1}
  Assume $G\in \mathcal{M}$ has a stable fixed point (say at zero)
  where it is analytic.  Then the difference equation (\ref{eq:01})
  has the Painlev\'e property {\bf iff} for some $a,b\in\CC$ with $|a|<1$,
  \begin{equation}
    \label{eq:concl}
    G(z)=\frac{az}{1+bz}
  \end{equation}
\end{theorem}
\z The proof is given in \S\ref{PT1}.
\begin{Remark}
  The Painlev\'e property {\em is not sensitive to which attracting
    fixed point of $G$ or its iterates is used in the analysis. This
  follows from the Proposition below.} 
\end{Remark}
Assume $p$ is another attracting fixed point of $G$ and let
$G_1(s)=G(p+s)-p$ ($G_1$ has an attracting fixed point at the origin).
\begin{Proposition}\label{Corol1}
  The difference equation (\ref{eq:01}) has the Painlev\'e property
  iff the difference equation $x_{n+1}=G_1(x_n)$ has the
  Painlev\'e property. Furthermore if $G$ has an iterate $G^{\circ m}$
  with an attracting fixed point where the conjugation map extends
  analytically to $\CC$ except for isolated singularities, then the
  same is true for any attracting fixed point of any iterate $G^{\circ
    k}$.

\end{Proposition}
\z This is shown in \S\ref{SProof112}.
\subsection{Failure of integrability test and barriers of singularities}\label{S19}
Conserved quantities are naturally defined as functions $C(x;n)$ with
the property $$C(x_{n+1};n+1)=C(x_n;n)$$

We now look at cases without the Painlev\'e property,  when $G$ is
a polynomial map.  We arrive at the striking conclusion that these
equations are not solvable in terms of functions extendible to the
complex plane, or on Riemann surfaces. The conserved quantities will
typically develop singularity barriers.

We use, in the formulation of the following theorem, a number of
standard notions and results relevant to iterations of rational maps;
these are briefly reviewed in the Appendix, \S\ref{Iter}.

\smallskip

\begin{theorem}\label{Barrier}
  Assume $G$ is a nonlinear polynomial with an attracting fixed point at
  the origin.  Denote by $\mathcal{K}_p$ the maximal connected component
  of the origin in the Fatou set of $G$. (It follows that $\mathcal K_p$
  is an open, bounded, and simply connected set).

  Then the domain of analyticity of $Q$ (see (\ref{invar})) is
  $\mathcal{K}_p$, and $\partial \mathcal{K}_p$ is a singularity barrier
  of $Q$.
\end{theorem}

\smallskip

\z This theorem is proved in \S\ref{PrTB}.
\subsection{Example: the logistic map}
The discrete logistic map is defined by
\begin{equation}
  \label{eq:log01}
x_{n+1}=a x_n (1-x_n)  
\end{equation}
The following result was proved by the authors in \cite{MKP}.
\begin{Proposition}\label{P4}
  The recurrence (\ref{eq:log01}) has the Painlev\'e property in
  Definition~\ref{DPP} iff $a\in\{-2,0,2,4\}$ (in which cases it is
  explicitly solvable). If $a\notin\{-2,0,2,4\}$ then the conserved
  quantity has barriers of singularities.
\end{Proposition}

\subsection{Application to the study of fractal sets}
The techniques also provide detailed information on the Julia sets of
iterations of the interval. 
\begin{theorem}\label{T5}
  Consider the equation (\ref{eq:log01}) for $a\in(0,1/2)$. 
  
  (i) There is an analytic function $G$, satisfying the functional
  relation
  \begin{equation}
    \label{frel}
    G(z)^2=aG(z^2)(1+G(z))
  \end{equation}
 which is a conformal map of
  the open unit disk $S_1$ onto $\{x^{-1}:x\in\text{ext}(\mathcal J)\}$
  where $\mathcal J $ is  the
  Julia set of (\ref{eq:log01}).
  
  (ii) $G$ is Lipschitz continuous of exponent $\log_2(2-a)$ in
  $\overline{S_1}$(the Lipschitz constant can be determined from the
  proof).
  
  (iii) $\partial S_1$ is a barrier of singularities of $G$. Near
  $1\in\partial S_1$ we have
\begin{equation}
  \label{formK}
  G(z)=\Phi\Big(\tau\Psi(\ln \tau)\Big)
\end{equation}
where
\begin{eqnarray}
  \label{defqt}
& \tau=\tau(z)=\ln(z^{-1})^{\log_2(2-a)}\nonumber \\
&\text{\rm $\Phi$ is
 analytic at zero, $\displaystyle \Phi(0)=\frac{a}{1-a}$,$\Phi'(0)=1$ }\\
&\text{\rm $\Psi$ is real analytic and periodic of period
 $\ln 2$.}\nonumber 
\end{eqnarray}
With $t=1-z$ we have
 \begin{equation}
    \label{pseudotrans}
      G=\frac{a}{1-a}+\sum_{l\in\ZZ}\sum_{k,m\in\NN}C_{l;k,m}
  t^{2\pi il\log_2(2-a)/\ln 2 +k\log_2(2-a)+m}
  \end{equation}
  where the series converges (rapidly) if $t$ and $|\arg t|$ are small.
\end{theorem}
\z This theorem is proved in \S\ref{S4}.
\begin{Note}
  The proof of Proposition \ref{P37} shows that the Lipschitz exponent
  is optimal. The theorem is valid for any $a<1$, and the proof is
  similar.
\end{Note}
\begin{Note}
  It follows from Theorem~\ref{T5} (iii) and (\ref{frel}) that every
  binary rational is a cusp of $\mathcal J$ of angle
  $\pi\log_2(2-a)$, see also Fig. 1.
\end{Note}

\figure {\label{Fig.1}}
\ifx\pdftexversion\undefined
$ $\hskip 4cm  \scalebox{0.5}[0.5]{\includegraphics[700,500][600,0]{050}}  \scalebox{0.5}[0.5]{\includegraphics[200,650][300,0]{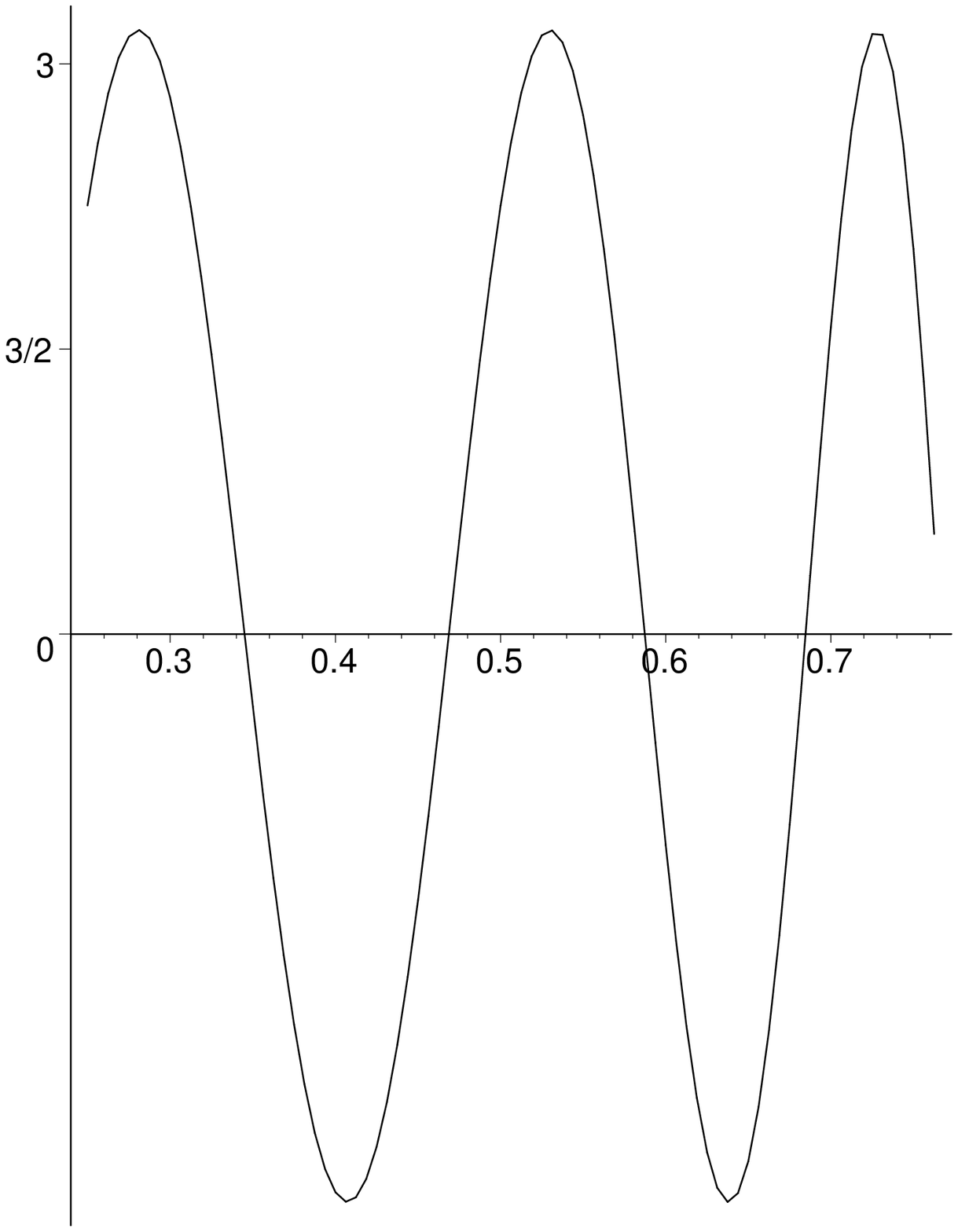}}\vskip 9cm
\else
 $ $\hskip -5cm\scalebox{0.5}[0.5]{\includegraphics{050}}
\vskip -10cm   \scalebox{0.5}[0.5]{$ $\hskip 9cm \includegraphics{Psi}}
\fi
\caption{(a) Julia set  for $G=\frac{1}{2}x(1-x)$.  The set $\mathcal{K}_p$ is the interior of the curve. (b) The function $10^9(\Psi(\ln\ln z_0)+c)$ for $a=\frac{1}{2},\,c=.079324389476$
(the plot relies on (\ref{limit}),  $N=300$).}
\endfigure

\section{General remarks on integrability}\label{hist1}

This problem has a long history, and the task of finding of {\em
  differential} equations solvable in terms of known functions was
addressed as early as the works of Leibniz, Riccati, Bernoulli, Euler,
Laplace, and Lagrange.  ``In the 18th century, Euler was defining a
function as arising from the application of finitely or infinitely
many algebraic operations (addition, multiplication, raising to
integer or fractional powers, positive or negative) or analytic
operations (differentiation, integration), in one or more variables''
\cite{BorelM}.  It was later found that some linear equations have
solutions which, although not explicit by this standard, have ``good''
{\em global} properties and can be thought of as defining new
functions.  To address the question whether nonlinear equations can
define new functions, Fuchs had the idea that a crucial feature now
known as the {\em Painlev\'e property} (PP) is the absence of movable
(meaning their position is solution-dependent, cf. \S\ref{Int})
essential singularities, primarily branch-points, see \cite{Fuchs}.
First order equations were classified with respect to the PP by Fuchs,
Briot and Bouquet, and Painlev\'e by 1888, and it was concluded that
they give rise to no new functions.  Painlev\'e took this analysis to
second order, looking for all equations of the form $u''=F(u',u,z)$,
with $F$ rational in $u'$, algebraic in $u$, and analytic in $z$,
having the PP \cite{Painleve1,Painleve2}. His analysis, revised and
completed by Gambier and Fuchs, found some fifty types with this
property and succeeded to solve all but six of them in terms of
previously known functions. The remaining six types are now known as
the Painlev\'e equations, and their solutions, called the Painlev\'e
transcendents, play a fundamental role in many areas of pure and
applied mathematics.  Beginning in the 1980's, almost a century after
their discovery, these problems were solved, using their striking
relation to linear problems\footnote{Some linear problems conducive to
  Painlev\'e equations were known already at the beginning of last
  century. In 1905 Fuchs found a linear isomonodromic problem leading
  to P$_{\rm VI}$.}, by various methods including the powerful
techniques of isomonodromic deformation and reduction to
Riemann-Hilbert problems \cite{Deift}, \cite{FlN}, \cite{Its}.

Sophie Kovalevskaya searched for cases of the spinning top having the
PP. She found a previously unknown integrable case and solved it in
terms of hyperelliptic functions. Her work \cite{Koval}, \cite{Koval2}
was so outstanding that not only did she receive the 1886 Bordin Prize
of the Paris Academy of Sciences, but the associated financial award
was almost doubled.

The method pioneered by Kovalevskaya to identify integrable equations
using the Painlev\'e property is now known as the {\em Painlev\'e test}.
Part of the power of the Painlev\'e test stems from the remarkable
phenomenon that equations passing it can generally be solved by some
method.  This phenomenon is not completely understood. At an intuitive
level, however, if for example all solutions of an equation are
meromorphic, then by solving the equation ``backwards,'' these solutions
and their derivatives can be written in terms of the initial
conditions. This gives rise to sufficiently many integrals of motion
with good regularity properties {\em globally} in the complex plane.

The Painlev\'e test has some drawbacks, notably lack of invariance under
transformations.  To overcome them, \cite{Kr-Cl} introduced the
poly-Painlev\'e test.

\section{Proofs}\label{S2}

\subsection{Proof of Theorem~\ref{T0}}\label{Uniq}

\subsubsection{Outline} The  idea of the proof is to use the convergence of (\ref{eq:transsS})
and its asymptotic properties to show that all terms $\mathbf
x_\mathbf k$ vanish.

\subsubsection{} We start with some preparatory results.

\begin{Remark}\label{R5}
  If $\mathbf x_\bfk\not\equiv 0$ then also $\bf X_k\not\equiv 0$ (see
  (\ref{eq:Bory})) so for small $p$ we have $\mathbf X_\mathbf k
  =\sum_{j=L_\bfk}^\infty \mathbf c_j p^j$ with $\mathbf c_{L_\bfk}\ne
  0$ for some $L_\bfk\ge 0$. By Watson's Lemma~\cite{Orszag}, for large
  $z$ in the right half plane we have
\begin{equation}
  \label{eq:Wat3}
  \mathbf x_\bfk\sim\sum_{j=L_\bfk}^\infty  \frac{\mathbf c_j j!} {z^{j+1}};\ \
  (\mathbf c_{L_\bfk}\ne 0)
\end{equation}
\end{Remark}

\begin{Remark}\label{finite}
  Since $\Re(\mu_i)>0$ we have $\Re(\boldsymbol{\mu}\cdot\bfk)
\rightarrow \infty$ as $\bfk \rightarrow \infty$.  Therefore for any
$K$, the sets of the form

\begin{equation}
  \label{eq:u3}
  \left\{\bfk\in\NN^{m_1}:\Re(\boldsymbol{\mu}\cdot\bfk)<K\right\};\ \ \left\{\bfk\in\NN^{m_1}:\Re(\boldsymbol{\mu}\cdot\bfk)=K\right\}
\end{equation}
are finite.
\end{Remark}
 
\begin{Definition}\label{D20}
  We define $S=\{\mathbf{k}:\mathbf X_\mathbf{k}\not\equiv 0\}$. We
  define inductively the finite  sets $T_i$ (cf. Remark~\ref{finite})
and  the numbers $M_i$ as follows:
  \begin{align}
    &T_0=\left\{\bfk\in S: \Re(\boldsymbol{\mu}\cdot\bfk)=\min_{\bfk\in
        S}\Re(\boldsymbol{\mu}\cdot\bfk)=:M_0\right\}\nonumber\\
   & T_1=\left\{\bfk\in S\setminus T_0:
      \Re(\boldsymbol{\mu}\cdot\bfk)=\min_{\bfk\in S \setminus
        T_0}\Re(\boldsymbol{\mu}\cdot\bfk)=:M_1\right\}\nonumber\\
     & \cdots\nonumber\\
   &  T_j=\left\{\bfk\in S\setminus T_0\ldots\setminus T_{j-1}:
      \Re(\boldsymbol{\mu}\cdot\bfk)=\min_{\bfk\in S\setminus
        T_0\ldots\setminus
        T_{j-1}}\Re(\boldsymbol{\mu}\cdot\bfk)=:M_j\right\}
\nonumber\\
 & \cdots
  \end{align}

    \z  Let also
    \begin{equation}
      \label{eq:defr}
      r_j=\max_{\bfk\in T_{j}}
      \Re(\mathbf a\cdot\bfk)
    \end{equation}
Note also that for some $\alpha>0$ we have
 \begin{equation}
      \label{eq:defr1}
      r_j\le \alpha M_j
    \end{equation}
\z Applying Remark~\ref{finite} again we see that
\begin{equation}
  \label{eq:union}
  \bigcup_{j=0}^{\infty}T_j=S
\end{equation}
 \end{Definition}

\begin{Lemma}\label{rearr}
  We have (see
  (\ref{eq:transsS})),
\begin{equation}
  \label{eq:decomp}
  \mathbf{x}(z)=\sum_{\bfk\in T_0}{\bf C}^{\bf k} e^{-{\bf k}\cdot
  \boldsymbol{\mu} z}
z^{{\bf k}\cdot {\bf a}}{\bf x}_\bfk(z)+O(e^{-M_1z}z^{r_1})\ \ \ \ 
(z\rightarrow +\infty)
\end{equation}
 \end{Lemma}
\begin{proof}
  We write 
  \begin{equation}
    \label{eq:dec2}
    \mathbf{x}(z)=\sum_{\bfk\in T_0}{\bf C}^{\bf k} e^{-{\bf k}\cdot
  \boldsymbol{\mu} z}
z^{{\bf k}\cdot {\bf a}}{\bf x}_\bfk(z)+\sum_{\bfk\in S\setminus T_0}{\bf C}^{\bf k} e^{-{\bf k}\cdot
  \boldsymbol{\mu} z}
z^{{\bf k}\cdot {\bf a}}{\bf x}_\bfk(z)
  \end{equation}
  The second series is uniformly and absolutely convergent for large
  enough $z\in\RR^+\!$ since it is bounded by the sub-sum of a
  (derivative of) a multi-geometric series
    \begin{equation}
      \label{eq:geoms}
      \sum_{\bfk\in S\setminus T_0}|{\bf A^{\bf k}}{\bf C}^{\bf k}z^{{\bf k}\cdot {\bf
      a}}| e^{-{\bf k}\cdot
  \Re(\boldsymbol{\mu}) z}
    \end{equation}
    Since (\ref{eq:geoms}) is absolutely
    convergent it can be thus be convergently rearranged as
 \begin{equation}
      \label{eq:geoms1}
      \sum_{j=1}^\infty e^{-M_jz} \sum_{\bfk\in  T_j}{|\bf A^{\bf
          k}}{\bf C}^{\bf k} z^{{\bf k}\cdot {\bf
      a}}| =\sum_{j=1}^\infty e^{-M_j z} z^{r_j} D_j(z)
    \end{equation}
    (see again Definition~\ref{D20} and Remark~\ref{finite}).  It is
    easy to see that $D_j(z)$ are nonincreasing in $z\in\RR^+\!$ and for large
    enough $z>0$ all products $z^{r_j}e^{-M_j z}$ are decreasing (cf. also
    (\ref{eq:defr1})). Therefore the convergent series
 \begin{equation}
   \label{eq:defs3}
\sum_{j=1}^\infty e^{-(M_j-M_1) z} z^{r_j-r_1} D_j(z)
 \end{equation}
is decreasing in $z>0$ and so
$$\sum_{j=1}^\infty e^{-M_j z} z^{r_j} D_j(z)\le Const. e^{-M_1 z}z^{r_1}$$
\end{proof}

\z {\bf Note}. A similar strategy could be also be used to show the
classical Weierstrass preparation theorem.

\subsubsection{}\label{Asy11} Assume first, to get a contradiction, that we have $\mathbf
x_0\not\equiv 0$ and so $\mathbf X_0\not\equiv 0$ so for small $p$ we
have $\mathbf X_\mathbf k =\sum_{j=m_0}^\infty \mathbf c_j p^j$ with
$\mathbf c_{m_0}\ne 0$. Then, since

$$\mathbf x(n)=\mathbf x_0(n)+O(e^{-M_1 n} n^{r_1})$$
\z and  by Remark~\ref{R5}
\begin{equation}
  \label{eq:Wat4}
  \lim_{n\rightarrow\infty} n^{-m_0-1}\mathbf x_0= (m_0+1)!c_{m_0}\ne 0
\end{equation}
 which contradicts $\mathbf x(n)=0$ for $n\in\NN$.
\subsubsection{}  Let now
\begin{equation}
  \label{eq:defR0}
R_0=\max\left\{\Re(\bfk\cdot\mathbf a-L_\bfk-1):\bfk\in T_0\right\}  
\end{equation}
and
\begin{equation}
  \label{eq:defTp0}
T'_0=\left\{\bfk\in T_0:\Re(\bfk\cdot\mathbf a-L_\bfk-1)=R_0\right\}  
\end{equation}
\begin{Lemma}\label{L18}
  We have
  \begin{equation} \label{eq:asympt} \mathbf
    x(z)=\sum_{\bfk\in T'_0}{\bf C}^{\bf k}
    c_{L_\bfk}L_\bfk!z^{\bfk\cdot\mathbf a-L_\bfk-1}e^{-\bfk\cdot\boldsymbol\mu
    z}+o\left(z^{R_0}e^{-M_0 z}\right) \ \ \mbox{for}\
    (z\rightarrow+\infty) \end{equation}
\end{Lemma}
\begin{proof}
  This is an immediate consequence of Remark~\ref{R5},
  Lemma~\ref{rearr}, 
and (\ref{eq:defR0}) and (\ref{eq:defTp0}).
\end{proof}

\subsubsection{Completion of the proof of Theorem~\ref{T0}} The
proof now follows, by reductio ad impossibile, from (\ref{eq:asympt}),
the assumption that $\mathbf x(n)=0$ for all large enough $n\in \NN$,
the fact that by construction all $c_{L_\bfk}$ are nonzero and the
following Lemma.

\begin{Lemma} \label{L20}

Let $d_{\bf k}\in\CC$. Then 
$$
\sum_{\bfk\in T'_0} d_\bfk n^{\bfk\cdot\mathbf
  a-M_\bfk}e^{-\bfk\cdot\boldsymbol\mu n}=o\left(n^{R_0}e^{-K_1
    n}\right)\ \ \  (\mbox{as }n\rightarrow\infty, n\in\NN)
$$
 iff all $d_\bfk$ are zero.
\end{Lemma}
\begin{proof}
  We now take $n_0=\mbox{card}(T'_0)$, $n$ large enough and note
  that $(n+j)^b=n^b(1+o(n^{-1}))$ if $j\le n_0$. Then a simple estimate shows
  that to prove the Lemma it suffices to show that the following
  equation cannot hold for all $0\le l\le n_0-1$ 
\begin{equation}\label{111}
\sum_{\bfk\in T'_0} d_\bfk e^{- (n+l)\bfk\cdot\boldsymbol\mu}=q_l
\end{equation}
 where
\begin{equation}
  \label{eq:defql}
  q_l=o(e^{-nM_0})\ \ \  (\mbox{as }n\rightarrow\infty, n\in\NN)
\end{equation}
\z If $n_0=1$ this is immediate. Otherwise, we may think of (\ref{111})
for $0\le l\le n_0-1$ as a system of equations for the $d_\bfk $ with
${\bf k}\in T'_0$. The determinant $\Delta$ of the system is a number of
absolute value $e^{-nlM_0}$ times the Vandermonde determinant of the
quantities $\{e^{- \bfk\cdot\boldsymbol\mu}\}_{\bfk\in T'_0}$. In
particular, for some $C>0$ independent of $n$ we have that
$e^{-nlM_0}|\Delta|$ is independent of $n$,
\begin{equation}
  \label{eq:as6}
e^{-nlM_0}|\Delta| =C\left|
  \prod_{\bfk_1\ne\bfk_2\in T'_0}(e^{-(\bfk_1-\bfk_2)\cdot\boldsymbol\mu}-1)\right|
\end{equation}
 and nonzero by (\ref{eq:nonres}).  Similarly, the minor $\Delta_\bfk$
of any $d_\bfk$ is bounded by $ D_\bfk e^{-n(l-1)M_0}$ with $D_\bfk $
independent of $n$. We get $d_\bfk=o(1)$ for large $n$ for all $\bfk\in
T'_0$, and so $d_\bfk=0$.
\end{proof}
\subsection{Remarks on first order equations}

\smallskip It turns out \cite{Levy} that for {\em first order
  autonomous} equations near an attracting fixed point, the series
$\tilde{x}_k$ of (\ref{eq:transszf}) are mere constants and the
transseries (\ref{eq:transszf}) are {\em classically convergent} for
large enough $n$ to actual solutions of the equation. This is a
consequence of the Poincar\'e equivalence theorem, see \cite{Levy}.

\begin{Note}
  If $|a|=1$ factorially divergent series do occur. In \S\ref{S31} we
  show how to use Borel summation instead of usual convergence when
  $a=1$.
\end{Note}

Assume for now that in (\ref{eq:01}) $G\in\mathcal{M}$ is analytic at
zero, $F(0)=F'(0)=0$ and $0<|a|<1$.  As we mentioned, there is a
one-parameter family of solutions presented as simple transseries of
the form
\begin{equation}
  \label{eq:2}
 x_n=x_n(C)=\sum_{k=1}^\infty e^{nk\ln a} C^k D_k
\end{equation}
with $D_k$ independent of $C$, which converge for large $n$. By
definition their continuation to complex $n$ is
\begin{equation}
  \label{eq:2z} x(z)=x(z;C)=\sum_{k=1}^\infty e^{z k\ln a} C^k D_k,
\end{equation}
 which is analytic for large enough $z$. To test for the Painlev\'e property, we
proceed to find the properties of $x(z)$ for those values of $z$ where
(\ref{eq:2z}) is no longer convergent, and then find the singular points
of $x(z)$.

\z {\bf Note.} In general, although (\ref{eq:2z}) represents a
continuous one-parameter family of solutions, there may be more
solutions.  We also examine this issue.

\subsubsection{Relation to properties of the conjugation map}
We can alternatively, and it turns out equivalently, define a
continuation as follows. By the Poincar\'e theorem \cite{Arnold} p.
99 there exists a unique map $\phi$ with the properties
\begin{equation}
  \label{assum1}
  \phi(0)=0,\ \  \phi'(0)=1\ \ \mbox{and} \ \phi \ \mbox{analytic at }\  0
\end{equation}
 and such that
\begin{equation}
  \label{eq:conjug}
  \phi(az)=G(\phi(z))=a\phi(z)+F(\phi(z))
\end{equation}
 The map $\phi$ is a {\em conjugation map} between (\ref{eq:01}) and
its linearization
\begin{equation}
  \label{Gens0}
  X_{n+1}=aX_n
\end{equation}
since, in view of (\ref{eq:conjug}),
\begin{equation}
  \label{Gens}
  x_n=\phi(Ca^n)
\end{equation}
for given $C$ and $n$ large enough, $x_n$ is a solution of the
recurrence (\ref{eq:01}).

We obtain a continuation of $x$ from $\NN$ to $\CC$ through
\begin{equation}
  \label{eq:con1}
x(z)=\phi(Ca^z)
\end{equation}
\begin{Lemma}\label{R1}
 (i) 
 For equations of type (\ref{eq:01}), the continuations (\ref{eq:2z})
and (\ref{eq:con1}) agree. 

\z (ii) $x(z;C)$ defined by (\ref{eq:2z}) has only isolated movable
singularities iff $\phi$ has only isolated singularities in $\CC$.
\end{Lemma}

\smallskip

\begin{proof}
  Indeed, $\phi$ is analytic at the origin, and a power series expansion
  for large $n$ of $\phi(Ca^n)$ leads to a solution of the form
  (\ref{eq:2}), which obviously solves (\ref{eq:01}). If $n_0$ is large
  enough, it is clear that (\ref{eq:2}) can be inverted for $C$ in terms
  of $x_{n_0}$ and we can also find $C'$ so that
  $x_{n_0}=\phi(C'a^{n_0})$.  On the other hand $x_{n_0}$ uniquely
  determines all $x_n$ with $n>n_0$. For equations of type
  (\ref{eq:01}), writing $x(z)=\phi(Ca^z)$ is thus tantamount to making
  the substitution $n=z$ in (\ref{eq:2}). Note that, $a^z$ is entire and
  $\phi$ is analytic at zero, and the presence of a singularity of
  $\phi$ which is not isolated is equivalent to the presence of a
  similar but {\em movable} singularity of $x(z)=\phi(Ca^z)$ since its
  position depends on $C$.
 \end{proof}

\subsubsection{Conserved quantities} The connection between $C$ and the equivalence map is seen as follows.  Near
an attracting fixed point, say $0$, we have a continuous one-parameter
family of solutions of (\ref{eq:01}) in the form (\ref{Gens}).

On the other hand the conjugation map $\phi$ is
invertible for small argument by (\ref{assum1}). We may then write
\begin{equation}
  \label{invar}
C=C(n,x_n)=\phi^{-1}(x_n)a^{-n}=:Q(x_n)a^{-n}
\end{equation}
where we see that $C(n,x_n)$ is a conserved quantity of (\ref{eq:01}),
and $Q=\phi^{-1}$ is analytic near zero.  Clearly any equation near a
stable fixed point is, in the sense of (\ref{invar}), {\em locally}
solvable.  Definition~\ref{DPP} requires however global properties.

Note first that, from the properties of $\phi$ (or from the constancy of
$C$), $Q$ satisfies the functional equation
\begin{equation}
  \label{eq:Q}
  Q(z)=a^{-1}Q(G(z))
\end{equation}

\subsection{End of proof of Lemma~\ref{L22} (a)}\label{PL4}
\begin{proof}
  For $i=1,2$, let $G_i\in M$, analytic in $\CC\setminus E_i$ and let
  $\CC\setminus E$ be the set of analyticity of $G_1\circ G_2$.  Then
  $E\subset\tilde{E}:=E_2\cup G_2^{-1}(E_1)$ is closed since the set
  of analyticity of any analytic function is open.  It remains to show
  $E$ is countable.  Since $G_2$ is not identically constant, for
  $x\not\in E_2$ there is a least $k=k(x)$ such that $G_2^{(k)}(x)\ne
  0$ and then $G_2$ has multiplicity exactly $k$ in a small disk $D_x$
  around $x$. Then $G_2^{-1}(E_1)\cap D_x$ is countable.  Since for
  every $x$ there is an open set $D_x$ such that $\tilde{E}\cap D_x$
  is countable it follows that $\tilde{E}$, thus $E$, is also
  countable. In the same way, for any $a\notin E_i$ we have that
  $G_i^{-1}(a)$ is countable.   For division, note that $1/G$ is
  defined wherever $G$ is defined and nonzero. Since $G$ is not a
  constant the same argument as above shows that $G^{-1}(0)$ is
  countable. \end{proof}\subsection{Proof of Theorem~\ref{T1}}\label{PT1}
\subsubsection{Notations}  In the following we will write $D_r(z_0)$ 
for the disk $\{z\in\CC:|z-z_0|<r\}$, $D_r$ will denote $D_r(0)$,
$\CC_\infty=\CC\cup\{\infty\}$. 

A number of notations, definitions and results in iterations of
rational maps used in the proof are reviewed in \S\ref{Iter}.

\begin{Proposition}\label{infmany}
  Let $R$ be a rational function of degree $d\ge 2$. Then $R$ has
  infinitely many distinct periodic points.
\end{Proposition}

\begin{proof}
  By definition, points of different period are distinct and by
  Lemma~\ref{classif} there are periodic points for every $n\ge 4$.
\end{proof}

\subsubsection{The ``if''  part of Theorem~\ref{T1}}
\label{SProof1}  In this direction 
the proof is trivial. Indeed, if $G$ is linear fractional, then the
general nonidentically zero solution of the equation (\ref{eq:01}) can
be obtained by substituting $x=1/y$ in (\ref{eq:01}) which then becomes
linear. We get
$$x_n=\left(Ca^{-n}+\frac{b}{a-1}\right)^{-1}$$
with the continuation
$x(z)=\left(Ce^{-z\ln a}+(a-1)^{-1}b\right)^{-1}$, a meromorphic
function.

\smallskip

\subsubsection{The ``only if''  part of Theorem~\ref{T1}}\label{SProof11} 
For the proof we will show that if $f$ has only isolated singularities,
and $f(az)=G(f(z))$, then $f$ itself is linear-fractional. Then $G$ is
also linear-fractional since $G(w)=f(a f^{-1}(w))$.
\begin{Lemma}\label{L21}
If $f$ has only isolated singularities and $f$ is not linear-fractional
then for any large enough $w$, the equation $f(z)=w$ has at least
two distinct roots.  
\end{Lemma}
\begin{proof}
  If $f$ is rational, then the property is immediate. Then assume that
  that $f$ is not rational, thus $f$ has at least one essential
  singularity, possibly at infinity \cite{Gonzales}. If $f$ has an
  essential singularity in $\CC$, then it is isolated by hypothesis
  and then the property follows from Theorem~\ref{BigPicard}. Then
  assume that $f$ has no essential singularity in $\CC$, thus infinity
  is the only essential singularity of $f$. If it is isolated then
  Theorem~\ref{BigPicard} applies again. Otherwise $f$ has infinitely
  many poles accumulating at infinity.  Since $f$ maps a neighborhood
  of every pole into a full neighborhood of infinity, any sufficiently
  large value of $f$ has multiplicity larger than one.
\end{proof}
 {\em End of proof of Theorem~\ref{T1}}.  Let now $G^{\circ m}$ be defined on
$\CC\setminus E_m$ and let $E=\cup_{m=1}^{\infty}E_m$; then $E$ is
countable and $G^{\circ m}$ is defined on $\CC\setminus E$ for any $m$. 

 Assume $f$ is not
linear-fractional and has only isolated singularities. We
let $z_1$ and $z_2$ be in $\CC\setminus E$ and such that
$f(z_1)=f(z_2)$, cf. Lemma~\ref{L21}.  Then
$f(az_1)=G(f(z_1))=G(f(z_2))=f(a z_2)$ and in general $f(a^n z_1)=f(a^n
z_2)$. But since $a^n z_1\rightarrow 0$ this contradicts (\ref{assum1}).

\subsection{Proof of Theorem~\ref{Barrier}}\label{PrTB}
  \begin{proof} 
    The fact that $\mathcal K_p$ is bounded for a nonlinear polynomial
    map follows from the fact that after the substitution $x=1/y$, the
    map $y_{n+1}=1/G(1/y_n)$ is attracting at $y=0$. Thus, cf.
    \cite{Beardon} Theorem 5.2.3 p. 83, $\mathcal{K}_p$ is simply
    connected. Let $a_1\in(|a|,1)$ and {\em let $D_\epsilon$ be a disk
      such that $|G(z)|<a_1 |z|$ for $z\in D_\epsilon$ and $Q$ is
      analytic in $D_\epsilon$}.

    By definition, for every $z_0\in\mathcal{K}_p$ there exists $m(z_0)$
    such that $G^{[m(z_0)]}(z_0)\in D_\epsilon$. Since $G^{[m(z_0)]}(z)$ is
    continuous in $z$, there is a disk $D_{\epsilon(z_0)}(z_0)$ such
    that $G^{[m(z_0)]}\left(D_{\epsilon(z_0)}(z_0)\right)\subset
    D_\epsilon$. It follows in particular that $\mathcal K_p$ is open.

    Since $\mathcal K_p$ is open and connected, it is arcwise connected.
    Let $z_0$ be arbitrary in $\mathcal{K}_p$ and let $C$ be an arc
    connecting $z_0$ to $z=0$. Since $C$ is compact and
    $$C\subset\bigcup_{z\in C}D_{\epsilon(z)}(z)$$ 

\z there is a finite
    subcovering 
        $$C\subset\mathcal{O}_C=\bigcup_{i=1}^ND_{\epsilon(z_i)}(z_i)$$
 with
    $z_i\in C$. Let $M$ be the largest of the $m(z_i), i=1,...,N$. Then,
    by construction, 
    \begin{equation}
      \label{eq:rezit}
G^{[M]}\left(\mathcal{O}_C \right)\in D_\epsilon       
    \end{equation}
 We see from (\ref{eq:Q}) that $aQ(z)=Q(G(z))=a^{-1}Q(G(G(z)))$ and
    in general, for $n\in\NN$,
 \begin{equation}
     \label{eq:period}
Q(z)=a^{-n}Q(G^{[n]}(z))     
   \end{equation}
   We define $Q(z)$ in $\mathcal{O}_C$ by $Q(z)= a^{-M}Q(G^{[M]}(z))$.
   By (\ref{eq:rezit}), and because (\ref{eq:period}) holds in
   $\mathcal{D}_\epsilon$, this unambiguously defines an analytic
   continuation of $Q$ from $D_\epsilon$ to $D_\epsilon\cup \mathcal
   O_C$. Since $\mathcal K_p$ is open and simply connected and since
   $Q$ is analytic near zero and can be continued analytically along
   any arc in $\mathcal K_p$, standard complex analytic results show
   that $Q$ is (single valued and) analytic in $\mathcal K_p$.

   For the last part, note that the boundary of $\mathcal{K}_p$ lies in
   the Julia set $J$, which is the closure of repelling periodic points
   (see Appendix, Lemma~\ref{repepo}).  Assume that $x_0$ is a repelling
   periodic point of $G$ of period $n$, and that $x_0$ is a point of
   analyticity of $Q$. Relation (\ref{eq:period}) implies that
   $Q(x_0)=0$ and that $Q'(x_0)=a^{-n}(G^{[n]})'(x_0)Q'(x_0)$ but since
   $|a|<1$ and $|{(G^{[n]}})'(x_0)|>1$ this implies $Q'(x_0)=0$.
   Inductively, in the same way we see that $Q^{(m)}(x_0)=0$ for all
   $m$, which under the assumption of analyticity entails $Q\equiv 0$
   which contradicts (\ref{assum1}).
   \end{proof}
\subsection{Borel summability of  formal invariant for logistic map  when $a=1$}
\label{S31}
We now consider an example which cannot be reduced to the previous
types, namely when $a=1$, and when therefore the Poincar\'e
equivalence theorem fails. In the recurrence
\begin{equation}
  \label{eq:log1}
  x_{n+1}=x_n(1-x_n)
\end{equation}
\z zero is a fixed point, and it can be shown in a rather
straightforward way that there are no attracting fixed points of this
map, or of any of its iterates. However, failure of the Painlev\'e property can be
checked straightforwardly, and Borel summability makes it possible to
analyze the properties of this equation rigorously.

A formal analysis of the Painlev\'e property is relatively
straightforward using methods similar to those in \cite{Invent}.  We
concentrate here on properties of the conserved quantities. 
 The
recurrence $a_{n+1}=a_n(1+a_n)^{-1}$ is exactly solvable and differs
from the logistic map by $O(a_n^3)$ for small $a_n$. The exact
solution is $n-a_n^{-1}=Const$, which suggests looking in the logistic
map case for a constant of the iteration in the form of an expansion
starting with $C=n-a_n^{-1}$.  This yields
\begin{equation}
  \label{eq:formali}
 C(n;v)\sim  n-{v}^{-1}-\ln v-\frac{1}{2}v-\frac{1}{3}\,{v}^{2}-{\frac {13}{36}}\,{v}^{3}-{
\frac {113}{240}}\,{v}^{4}+\cdots
\end{equation}
which is indeed a formal invariant, but the associated series is
factorially divergent as will appear clear shortly. Nevertheless we can
show that the expansion is Borel summable to an actual conserved
quantity in a sectorial neighborhood of $v=0$.

\begin{theorem}\label{T3} There is a conserved quantity $C$ defined near
the origin in $\CC\setminus\RR^-$, of the form
$$C(n;v)=n-v^{-1}-\ln(v)-R(v)$$
 where $R(v)$ has a Borel summable series at the origin in any
direction in the open right half plane. $R(v)$ has a singularity barrier
touching the origin tangentially along $\RR^-$. This singularity barrier
is exactly the boundary of the Leau domain of (\ref{eq:log1}).
\end{theorem}

 We let
\begin{equation}
  \label{eq:R(v)}
  C(n;v):=n-v^{-1}-\ln v-R(v)
\end{equation}
and impose the condition that $C$ is constant along
trajectories. This yields
\begin{equation}
  \label{eq:recR}
  R(v)=R(v-v^2)+ \frac{v}{1-v}+\ln(1-v)
\end{equation}
where the RHS of (\ref{eq:recR}) is $R(v-v^2)+O(v^2)$. The
substitution 
\begin{equation}
  \label{eq:defR}
R(v)=h(v^{-1}-2)  
\end{equation}
followed by $v=1/(x+1)$ yields
\begin{equation}
  \label{eq:eqh}
  h(x-1)=h\left(x+x^{-1}\right)+\frac{1}{x}+\ln\left(\frac{x}{x+1}\right)
\end{equation}
which by formal expansion in powers of $x^{-1}$  becomes
\begin{equation}
  \label{eq:hexp}
  h(x-1)=\sum_{k=0}^{\infty}\frac{h^{(k)}(x)}{k!}x^{-k}+\frac{1}{x}+\ln\left(\frac{x}{x+1}\right)
\end{equation}
\subsubsection{Proof of Theorem~\ref{T3}}\label{S32}

\begin{Proposition}
\newcounter{bb}
\label{Cor1} $ $
{
\begin{list}{\roman{bb})}{\usecounter{bb}\setlength{\leftmargin}{0pt}}
\item \newcounter{p}\setcounter{p}{\value{bb}} $R(v)=h(v^{-1}-2)$ has a
  Borel summable series at the origin along $\RR^+\!$. More precisely, $h(x)$ can be written in the form
  \begin{equation}
    \label{eq:reph}
    h(x)=\int_0^\infty e^{-px}H(p)dp
  \end{equation}
  and where $H(p)$ is analytic at zero and in the open right half plane
  $\mathbb{H}=\Re(p)>0$ and has at most exponential growth along any ray
  towards infinity in $\mathbb{H}$.

\item \newcounter{pi}\setcounter{pi}{\value{bb}} $h$ is analytic in a
  region of the form $\{x:\arg(x)\ne \pi: |x|\ge \nu(\arg(x))\}$. The
  function $\nu$ is continuous in $(-\pi,\pi)$. (The expression of
  $\nu:(-\pi,\pi)\mapsto\RR^+\!$ will follow from the proofs below.)

\item \newcounter{pii}\setcounter{pii}{\value{bb}} By (\ref{eq:reph})
and Watson's Lemma~\cite{Orszag}, $h$ has an asymptotic power series for
large $x$, $h(x)\sim\sum_{k=0}^{\infty}H^{(k)}(0)x^{-k}$, which is a
formal solution of (\ref{eq:hexp}).

\item \newcounter{piii}\setcounter{piii}{\value{bb}}\label{iii}The function $R(v)$
  is analytic in a region near the origin, the origin excluded, of the
  form $\mathcal{V}=\{v: \arg(v)\ne \pi, 0<|v|<\nu^{-1}(\arg(\phi)\}$.
   By (\roman{pii}) the relation (\ref{eq:formali})
  is an asymptotic expansion for small $v\in\mathcal{V}$, and from
  (\ref{eq:reph}) the power series contained there is Borel summable.

\item \newcounter{piv}\setcounter{piv}{\value{bb}} The function $R$
  given by (\ref{eq:defR}) satisfies (\ref{eq:recR}).

\item \newcounter{pv}\setcounter{piv}{\value{bb}} The function $R$ is
analytic in $L_f$, the Leau domain of $f$, and has a singularity barrier
on the Julia set of $f$.

\end{list}
}
\end{Proposition}

\begin{proof}

 The formal inverse Laplace transform of (\ref{eq:hexp}) is the
equation
\begin{equation}
  \label{eq:eqH}
  (e^p-1)H=\frac{1-e^{-p}-p}{p}+\sum_{k=1}^\infty \frac{(-p)^k}{k!}H*1^{*k}
\end{equation}
where * denotes the Laplace-type convolution
$$F*G=\int_0^p F(s)G(p-s)ds$$
and $F^{*k}$ is the convolution of $F$ with itself $k$ times.
We rewrite (\ref{eq:eqH}) in the form
\begin{equation}
  \label{eq:eqH1}
  H=\frac{1-e^{-p}-p}{p(e^p-1)}+\frac{1}{(e^p-1)}\sum_{k=1}^\infty
  \frac{(-p)^k}{k!}H*1^{*k}=H_0+\mathfrak{A} H
\end{equation}
where $\mathfrak{A}$ is a linear operator. We show now that this
equation is contractive in an appropriate space of functions. Let
$\nu>0$ and let $\mathcal{A}$ be the space of functions $F$ analytic in
a neighborhood $\mathcal{N}$ of $[0,\infty)$ in the complex plane, with
$F(0)=0$, in the norm $\|F\|_\nu:=\sup_{\mathcal{N}}|e^{-\nu
|p|}F(p)|$. We choose $a\in (0,2\pi)$, $\epsilon$ small and
\begin{equation}
  \label{eq:defN}
  \mathcal{N}=\Big\{p:|p|\le \epsilon\Big\}\cup \left\{p: \arg(p)\in
  \left(-\frac{\pi}{2}+\epsilon,
\frac{\pi}{2}-\epsilon\right)\right\}
\end{equation}
Since the norm $\|\cdot\|_\nu$ restricted to compact sets is equivalent
to the usual sup norm, it is easy to check that $\mathcal{A}$ is a
Banach space.
\begin{Proposition}
  For large enough $\nu$, the equation (\ref{eq:eqH1}) is contractive
in $\mathcal{A}$ in the norm $\|\cdot\|_\nu$.
\end{Proposition}
First, it is easy to see that $H_0\in\mathcal{A}$. If
$f\in\mathcal{A}$
then 
\begin{multline}
  \label{eq:Af}
  \sum_{k=1}^\infty
  \frac{(-p)^k}{k!}f*1^{*k}=
\sum_{k=1}^\infty\frac{(-p)^k}{k!}\int_0^p f(s)
\frac{(p-s)^{k-1}}{(k-1)!}ds\\=\sum_{k=1}^\infty\frac{(-1)^kp^{2k}}{k!(k-1)!}\int_0^1 f(pt)
(1-t)^{k-1}dt=\sum_{k=1}^\infty\frac{(-1)^kp^{2k}}{k!(k-1)!}\int_0^1 f(p(1-t))
t^{k-1}dt
\end{multline}
It is immediate that if $p$ is in a compact set $\mathcal{K}$ and $f$
is analytic in $\mathcal{K}$ then the sum in (\ref{eq:Af}) is uniformly
convergent in $\mathcal{K}$ and analytic in $p$. Furthermore the sum is
$O(p^3)$ for small $p$ since $f\in\mathcal{A}$. Now we see that
\begin{multline}
  \label{eq:est3}
  \left|e^{-\nu |p|}\sum_{k=1}^\infty\frac{(-1)^kp^{2k}}{k!(k-1)!}\int_0^1 f(p(1-t))
t^{k-1}dt\right|\\=\left|
\sum_{k=1}^\infty\frac{(-1)^kp^{2k}}{k!(k-1)!}\int_0^1 e^{-\nu |p|(1-t)}f(p(1-t))
t^{k-1}e^{-\nu |p|t}dt\right|\\
\le
\|f\|_\nu\sum_{k=1}^\infty\frac{|p|^{2k}}{k!(k-1)!}
\int_0^{1}t^{k-1}e^{-\nu |p|t}dt\le
\|f\|_\nu\sum_{k=1}^\infty\frac{|p|^{2k}}{k!(k-1)!}
\int_0^{\infty}t^{k-1}e^{-\nu |p|t}dt
\\ = \|f\|_\nu\sum_{k=1}^\infty\frac{|p|^{k}}{k!\nu^k}\le 
\|f\|_\nu\frac{|p|}{\nu}e^{|p|/\nu}
\end{multline}
and thus
\begin{equation}
  \label{eq:norm}
  \|\mathfrak{A}\|\le Const\,\nu^{-1}
\end{equation}
for sufficiently large $\nu$, where we took into account the
exponential decrease of $(e^p-1)^{-1}$ for large $p$ in $\mathcal{N}$.
Thus the equation has a unique fixed point  $H\in\mathcal{A}$. In
particular the Laplace
transform $h(x)=\mathcal{L}H=\int_0^\infty e^{-xp}H(p)dp$ is well defined and
analytic in the half-plane
$\Re(x)>\nu$. It is now immediate to check that $h(x)$ satisfies
the equation (\ref{eq:hexp}).

\section{Julia sets for the map (\ref{eq:log01}) for $a\in(0,1)$}\label{S4}
It is convenient to analyze the superattracting fixed point at infinity; the
substitution $x=1/y$ transforms (\ref{eq:log01}) into
\begin{equation}
  \label{logisty}
  y_{n+1}=-\frac{y_n^2}{a(1-y_n)}
\end{equation}
For small $y_0$, the leading order form of equation (\ref{logisty}) is
$y_{n+1}=-a^{-1}{y_{n}^2}$ whose solution is $-y_0^{2^n}a^{-2^{n-1}-1}$. It is
then convenient to seek solutions of (\ref{logisty}) in the form
$y_n=-G(y_0^{2^n}a^{-2^{n}})$ whence the initial condition implies $G(0)=0$,
$G'(0)=a$. Denoting $y_0^{2^n}a^{-2^{n}}=z$, the functional relation
satisfied by $G$ is
\begin{equation}
  \label{eq:112}
 G(z^2)=\frac{G(z)^2}{a(G(z)+1)}; \ \ \ \ G(0)=0,\ G'(0)=a
\end{equation}
\begin{Lemma}[\cite{MKP}]\label{L11}
There exists a unique analytic function $G$ in the neighborhood of the
origin satisfying (\ref{eq:112}). This $G$ has only isolated singularities in $\CC$ if and
only if $a\in\{-2,2,4\}$. In the latter case, 
(\ref{eq:log01}) can be solved explicitly. 

If $a\not\in\{-2,2,4\}$ then the unit disk is a barrier of singularities
of $G$.

\end{Lemma}
\begin{Lemma}\label{L14} 
  $G$ is analytic in the open unit disk $S_1$ and Lipschitz continuous in
  $\overline{S_1}$. 
\end{Lemma}
\begin{proof} Lemma \ref{L11} proved in \cite{MKP}
  guarantees the existence of some disk $S_r$ centered at zero, of
  radius $r\le 1$, where $G$ is analytic and it is shown that inside
  that disk we have (cf.  also \ref{eq:112})
\begin{equation}
  \label{113TMP}
  G(z)=U(G(z^2)); \ \ 2U(s):=s+(a^2s^2+4s)^{\frac{1}{2}}
\end{equation}
(with the choice of branch consistent with $G(0)=0, G'(0)=a$). If $r<1$
then (\ref{113TMP}) provides analytic continuation in a disk of radius
$r^{\frac{1}{2}}>r$ if $a^2 G(z)^2+4aG(z)\ne 0$ in $S_r$. 
\begin{Note}\label{N0} $G(z_0)=0$ in $S_1$ iff $z_0=0$.
\end{Note}
\z Indeed, assume $0\ne z_0\in S_r$ and $G(z_0)=0$. Then we find from
(\ref{eq:112}) that $G(z_0^{2^n})=0$ which is impossible since $G$ is analytic
at zero and $G'(0)=a$. $\Box$

We are left to examine the possibility
$G(z_0)=-4a^{-1}$ with $z_0\in S_r$.
\begin{Note}\label{N1} $R(x)=\frac{x^2}{a(x+1)}$ is well defined and increasing on the interval
  $\left(-\infty,-4a^{-1}\right)$.
\end{Note}
The assumption $G(z_0)=-4a^{-1}$ thus implies that the values
$G(z_0^{2^n})$ are in $\RR^-$ and decrease in $n$, again impossible
$G$ is analytic at $0$ and $G(0)=0$.

We now show $G$ is bounded in $S_1$. Indeed, by (\ref{113TMP}) we have
\begin{equation}
  \label{ineq}
|G(z)|\le U(|G(z^2)|)\end{equation}
on the other hand, a calculation shows that 
\begin{equation}
  \label{basicineq}
 U(s)\le\frac{a}{1-a}\ \ \text{for }s\in
\left[0,  \frac{a}{1-a}\right]
\end{equation}
Since $G(0)=0$ and $G$ is analytic in $S_1$, (\ref{ineq}) and
(\ref{basicineq}) imply that 
\begin{equation}
  \label{estimsup}
  \sup_{z\in S_1}|G(z)|\le \frac{a}{1-a}
\end{equation}
We next prove that $G$ is injective. As a first step we have the following: 
\begin{Note}\label{N2}  $G'\ne 0$ in
$S_1$.
\end{Note}
Indeed, otherwise differentiating (\ref{eq:112}) shows there would exist a
sequence $z_n\rightarrow 0$ such that $G'(z_n)=0$.$\Box$

Now, $G$ is injective in a neighborhood of the origin since $G'(0)=a$. Let
then $z_1\in S_1$ be a point of smallest modulus such that there exists
$z_2\ne z_1\in S_1$ with $G(z_1)=G(z_2)$. For $z_1$ to exist, we need, again
by (\ref{eq:112}) that $z_1^2=z_2^2$ and thus $z_1=-z_2$. Since $G'\ne 0$, by
the open mapping theorem, the image under $G$ of arbitrarily small disks
around $z_1$ and $-z_1$ overlap nontrivially. For some $C$ and any $\epsilon$
there exist therefore infinitely many $z_i$ with $|z_1-z_i|<\epsilon$ such
that $G(z_i)=G(z'_i)$ and $|z'_i-(-z_1)|<C\epsilon$. The same argument using
(\ref{eq:112}) shows that $z'_i=-z_i$. But since $G(z)=G(-z)$ for infinitely
many $z\in S_1$ accumulating at $z_1$, then $G$ would be even, which is not the
case since $G'(0)=a$.
We now need two lower bounds.
\begin{Proposition}
  \label{P35} For $a\in \left(0,\frac{1}{2}\right)$ 
$$(1-|z|)^{1-\log_2(2-a)} G'(z)$$ is bounded in $S_1$.
\end{Proposition}
\begin{proof}
  The function $H=1/G$ which, by Note~\ref{N0}, is analytic in $S_1\setminus 0$
satisfies 
\begin{equation}
  \label{eH3}
  H(z^2)=aH(z)(1+H(z))
\end{equation}
Let 
$$m_n=\max\{|H(z)|:|z|\in\left[2^{-\frac{1}{2^n}},
  2^{-\frac{1}{2^{n+1}}}\right]\}$$
 Eq. (\ref{eH3}) gives
$$m_{n+1}\le \frac{1}{2}+\sqrt{\frac{1}{4}+\frac{m_n}{a}}$$
and it easy to see that this implies
\begin{equation}
  \label{eq:e4}
  \limsup_{n\rightarrow\infty}m_n\le 1+a^{-1}
\end{equation}
 We have
\begin{equation}
  \label{2step}
G'(z)=2az\frac{G'(z^2)(1+G(z))^2}{G(z)(2+G(z))}  
\end{equation}
so that
 $$|G'(z)|\le |G'(z^2)|\max_{(1-a)|y|\le a}\left|\frac{2a(1+y)^2}{y(2+y)}\right|= \frac{2}{2-a}|G'(z^2)|$$
 if $a\le 1/2$ from which Proposition ~\ref{P35} follows
 immediately. 
\end{proof}
 \begin{Note}
   A straightforward way to extend the result for larger values of
   $a<1$ is to replace (\ref{2step}) by a corresponding equality
   obtained from a higher order  iterate of (\ref{eH3}).
 \end{Note}

\begin{Lemma}\label{L15} 
$G$ is gives a conformal transformation of $S_1$ onto a
  bounded region $\mathcal{K}_p$, whose boundary $\partial \mathcal{K}_p$ is
  a Lipschitz continuous nowhere differentiable curve. 
\end{Lemma}

\section{Behavior at the singularity barrier}
\begin{Proposition}
  \label{P37}
  There is $\delta>0$, a real analytic function $\Psi$, periodic of
  period $\ln 2$ and an analytic function $\Phi, \Phi'(0)=1$ such that
  for $|\arg(1-z)|<\delta$ (\ref{formK}) holds.
\end{Proposition}

\begin{proof}
  Let $\omega=2\pi/\ln 2$, $\beta=\log_2(2-a)$.  With $z_0\in (0,1)$
  and $z_n=z_0^{1/2^n}$, the sequence $G_n=G(z_n)$ is increasing and
  bounded by $L$, see (\ref{113TMP}). It follows immediately from
  (\ref{113TMP}) and (\ref{ineq}) that
\begin{equation}
  \label{conv3}
  L-G_n:=\delta_n\downarrow 0 \ \ \text{as}\ \ n\rightarrow\infty\ \ (L:=\frac{a}{1-a})
\end{equation}
From (\ref{eq:e4}) we have, with $C_1=\frac{(1-a)^3}{a(2-a)}+C_2$, $C_2=1-a$
\begin{equation}
  \label{itdel}
  \delta_{n+1}=\frac{1}{2-a}\left[\frac{1-C_2\delta_{n+1}}{1-C_1\delta_{n+1}}\right]\delta_n
\end{equation}
Eqs. (\ref{conv3}) and (\ref{itdel}) imply that for any
$\epsilon>0$ we have
\begin{equation}
  \label{geomd}
  \delta_n=o\left((2-a-\epsilon)^{-n}\right)\ \ \ \text{as } n\rightarrow\infty
\end{equation}
 Let 
 \begin{equation}
   \label{defdn}
   \delta_n=\ln^\beta(1/z_n)e^{\theta_n}=2^{-n\beta}
\ln^\beta(1/z_0) e^{\theta_n} 
 \end{equation}
cf. (\ref{conv3}). Now
$$\left|e^{\theta_n-\theta_{n+1}}-1\right|=\frac{(C_1-C_2)\delta_{n+1}}{1-C_2\delta_{n+1}}=O(\delta_n)\ \ \ \text{as } n\rightarrow\infty$$
\z and by (\ref{geomd}), $\theta_n$ is convergent, $\theta_n\rightarrow\Theta$. Since
$\theta_{n+1}-\theta_n\rightarrow 0$ it follows that
\begin{equation}
  \label{per1}
  \Theta(z_0^2)=\Theta(z_0)
\end{equation}
\subsubsection{Analyticity}\label{611}  We let $1-z_1$ be sufficiently small so that 
\begin{equation}
  \label{iter4}
  \delta_n\le c \alpha^n
\end{equation}
\z with $\alpha<1$ and $c$ small enough so that the term in square
brackets is sufficiently close to one for all $n\ge 0$ and
$|z_0-z_1|\le \epsilon_1$ (cf.  (\ref{geomd})), this amounts to a
shift in $n$). If $\epsilon_1$ is small enough, then it is easy to
check that equation (\ref{itdel}) is a contractive mapping in the in
the ball of radius $c$ $S_{\epsilon_1}=\{\zeta:|\zeta\le \epsilon_1\}$
in Banach space $l_{\infty,\alpha}(\NN) $ of vectors
$\mathbf{v}(n;\zeta)$ analytic in $\zeta=z_0-z_1$ with respect to the
norm
$$\|\mathbf{v}\|=\sup_{n\ge 1;|\zeta|\le \epsilon_1}|\mathbf{v(n,\zeta)\alpha^{-n}}|$$
and local analyticity in a neighborhood of the interval $[z_0, \sqrt{z_0}]$.
By periodicity, real analyticity follows immediately and relation (\ref{per1})
is preserved.

\end{proof}

\subsubsection{End of the proof of Theorem~\ref{T5} (iii)}
We use the information obtained in \S\ref{611}. Let  $e^{\theta_n}=(1+w_n)e^\Theta$; given $\delta>0$ we choose $n_0$ large enough
and $\epsilon_2$ so that $|w_n(z_0)|<\delta$ if
$|z-z_0|<\epsilon_2$ and $n\ge n_0$. We let ,  $h=e^{2\Theta}$, 
$\varepsilon_n=2^{n\beta}$, $s=\ln^\beta(1/z_0)$, $c=C_1-C_2, C=c-C_2$ and obtain
\begin{multline}
  \label{ctrw}
  w_n=\frac{Ce^{2\Theta}s\varepsilon_n}{1-\varepsilon_n C_2 e^{2\Theta}s}\\
+w_{n+1}\frac{1+2Ce^{2\Theta}s\varepsilon_n-C C_2e^{4\Theta}s^2\varepsilon_n^2+w_{n+1}Cse^{2\Theta}\varepsilon_n(1-\varepsilon_nC_2se^{2\Theta})}{1-2\varepsilon_n se^{2\Theta}+2\varepsilon_n^2 s^2e^{4\Theta} -w_{n+1}C_2se^{2\Theta}(1-\varepsilon_n C_2 se^{2\Theta})}
\end{multline}
As in \S\ref{611}, a contractive mapping argument shows that
$\mathbf{w}=(w_n,w_{n+1},...)$ is analytic in $se^{2\Theta}$, if $s$
is small enough. The conclusion now follows from the definition
$$G(z_0^{2^{-n_0}})=L+s\delta_{n_0}$$
and (\ref{defdn}), (\ref{per1}),
\S\ref{611} and the substitution
$e^{2\Theta(\cdot)}=\Psi(\ln(\ln(\cdot))$.
Formula (\ref{pseudotrans}) follows immediately from (\ref{formK}).
\end{proof}
\begin{Note}
  With $z_n=z_0^{1/2^n}$, $\tau_n=\tau(z_n)$ (cf. (\ref{formK}) and
  $g_n=G(z_n)-L$ we have
\begin{equation}
  \label{limit}
  \Psi(\ln\ln z_0)=\lim_{N\to \infty}\frac{\displaystyle \frac{g_{N+1}}{\tau_{N+1}}-\frac{g_{N}}{\tau_{N}}}{\tau_{N+1}-\tau_N}
\end{equation}

\end{Note}

%%      ---------------------------------------------------------------------
%%      ---------------------------- BODY OF PAPER --------------------------
%%      ---------------------------------------------------------------------

%%      Please input or insert the body of your paper here.

%%      ---------------------------------------------------------------------
%%      ------------------------- APPENDIX (OPTIONAL) -----------------------
%%      ---------------------------------------------------------------------

\appdx[Iterations of rational maps]
\label{Iter}
We introduce a number of definitions and results for iterations of
rational maps, which are treated in much more detail and generality in
\cite{Steinmetz} and \cite{Beardon}. We shall illustrate the main
concepts on the simple case $G=ax(1-x)$. In Figure 1, the interior (in
the complex plane) of the fractal curves is a set invariant under $G$
and with the further property that starting with $z_0$ inside the m-th
iterate of $G$ at $z_0$, $G^{\circ m}(z_0)$, converges to zero as
$m\rightarrow\infty$.  These are {\em stable fixed domains} of $G$.

Consider the polynomial map $G$.  A Fatou domain of $G$ is a stable
fixed domain $V$ of $G$ characterized by the property that $G^{\circ n}$
converges {\em in the chordal metric} on the Riemann sphere
$\mathbb{C}_\infty$ to a fixed point  of $G$, locally uniformly in
$V$.

\begin{Definition}[\cite{Beardon}, p. 50] Let $G$ be a non-constant rational
  function. 
The {\bf Fatou set} of $G$ is the maximal open subset of
$\mathbb{C}_\infty$
on which $\{G^{\circ n}\}$ is equicontinuous and the {\bf Julia set} of $G$ is
its complement in $\mathbb{C}_\infty$.
  
\end{Definition}

A Fatou domain is a {\bf Leau} domain (or a parabolic basin) if
$x_0\in\partial V$ and the multiplier of $x_0$ (the derivative at $x_0$)
is $\lambda=1$\footnote{\cite{Steinmetz}, p.  54}. In Figure 1 this
happens for $a=1$.

The Julia set can be characterized by the following property.

\begin{Lemma}[\cite{Beardon}, p. 148] Let $G$ be a rational map of degree $d$,
  (cf. Definition~\ref{defdeg}) where $d\ge 2$. Then $J$ is the derived
  set\footnote{By definition the derived set of a set $E$ consists
  exactly in the points $z$ which are limits of sequences $\{z_n\}$
  where the $z_n\in E$ are distinct.} of the periodic points of $G$.
\end{Lemma}

\z Under the assumptions above, we have 
\begin{Lemma}[\cite{Beardon}, p. 148]\label{repepo} $J$ is the closure of the
  repelling points of $G$.
\end{Lemma}

\begin{Definition}[\cite{Beardon}, p. 30.]\label{defdeg}
  If $R=P/Q$ where $P$ and $Q$ are polynomials, then the degree of the
  rational function $R$ is $\max\{\deg(P),\deg(Q)\}$.
\end{Definition}
\begin{Definition}[\cite{Beardon}]\label{defperi}
  If $R$ is a rational function and $R^{\circ m}$=$R\circ R\circ\cdots\circ
  R$ $n$ times, then a periodic point of period $n$ of $R$ is a point
  $z$ such that $R^{\circ m} z=z$ and $R^{\circ m} z\ne z$ if $m<n$. A periodic
  point of $R$ is a point of some period $n\ge 1$.

\end{Definition}

\z We also use the following result of I. N.  Baker:

\begin{Lemma}[\cite{Baker}, \cite{Beardon}]\label{classif}
  Let $R$ be a rational function of degree $d\ge 2$, and suppose that
  $R$ has no periodic points of period $n$. Then $(d,n)$
is one of the pairs

$$(2,2), (2,3), (3,2), (4,2)$$
\z (moreover, each such pair does arise from some $R$ in this way).
\end{Lemma}

\subsubsection{Further results used in the proofs}
\begin{theorem}[Big theorem of Picard, local formulation \cite{Rudin},
  \cite{Gonzales}]\label{BigPicard} If $f$ has an isolated singularity at a point $z_0$
  and if there exists some neighborhood of $z_0$ where $f$ omits two
  values, then $z_0$ is a removable singularity or a pole of $f$.
\end{theorem}

\begin{theorem}[Picard-Borel, \cite{Nevanlinna}]\label{PB}
  If $\phi$ is any nonconstant function meromorphic in $\CC$, then
  $\phi$ avoids at most two values (infinity included).
\end{theorem}
All we need in the present paper is that at most two finite values are
excluded. This is immediately reduced to the more familiar Picard theorem
by noting that if $\lambda$ is an excluded value of $f$ then $1/(f-\lambda)$
is entire.

\subsection{Proof of Proposition~\ref{Corol1}}\label{SProof112} By Theorem~\ref{T1},
(\ref{eq:01}) does not have the Painlev\'e property at some stable fixed point iff $G$
is not linear-fractional, in which case (\ref{eq:01}) fails to have the
Painlev\'e property at any other stable fixed point. More generally,
Proposition~\ref{Corol1} follows from the following result.

\begin{Lemma}\label{L12}
  If $G^{\circ m}$ is of the form (\ref{eq:concl}) then $G$ is of the form
(\ref{eq:concl}).
\end{Lemma}

\begin{proof}
  Since (\ref{eq:concl}) is one to one, the conclusion follows from the
  remark that if $G$ is not linear-fractional, then $G(z)$ has
  multiplicity greater than one for all sufficiently large $z$ (and then
  the same holds for $G^{\circ m}(z)$).  Indeed, assume that $G$ is not
  linear-fractional.  If $G$ is rational, then the conclusion is
  obvious. If the set of singularities of $G$ is finite, then they are
  all isolated and at least one is an essential singularity (otherwise
  $G$ is rational \cite{Gonzales}) and Theorem~\ref{BigPicard} applies.

So we may assume the set of singularities is infinite. Since by
  assumption this set is closed and countable, it contains infinitely
  many isolated points. (Indeed, a set which is closed and dense in
  itself, i.e. a {\em perfect set}, is either empty or else
  uncountable.) Then if $G$ has an isolated essential singularity,
  Theorem~\ref{BigPicard} applies, and if not then there are infinitely
  many poles of $G$. In the latter situation any sufficiently large
  value of $G$ has multiplicity larger than one since $G$ maps a
  neighborhood of every pole into a full neighborhood of infinity.
\end{proof}

\z {\em Completion of proof of Proposition~\ref{Cor1}}. Part
(\roman{pi}) merely follows from the formula (\ref{eq:reph}) and
elementary contour deformation in the integral. Parts (\roman{pii}) and
(\roman{piii}) are straightforward.

After the transformation $v=-u+1/2$ the iteration associated to our
map $f$ is equivalent to that of the quadratic map $q(u)= u^2+1/2$.

Part (\roman{piv}) follows  from the following Lemma.

\begin{Lemma}[\cite{Steinmetz}, p. 174] The Leau domain of $q$,   is the filled in
  (interior of the) Julia set $\mathcal{K}_p$ of $q$.
\end{Lemma}
 {\em Proof of Proposition~\ref{Cor1}(\roman{piv})}.  Let
   $H(v)=R(v)+v^{-1}+\ln v$, defined and analytic in $\mathcal{V}$. By
   definition we have $H(v_{n+1})=H(v_n)+1$ i.e.
\begin{equation}
  \label{eq:26}
 H(v)= H(f(v))-1=H(v-v^2)-1
\end{equation}
and clearly $R$ and $H$ have the same type of singularities in
$\CC\setminus\RR^-\setminus\{0\}$. 

If $z_0\in L_f$ we have by definition $|z_n|=|f^{\circ m} (z_0)|\rightarrow
0$.  Then, we choose $\epsilon$ small enough and $N$ so that
$|z_n|<\epsilon$ for $n>N$. Since we must have for some $n>N$ that
$|z_{n+1}|<|z_n|$, then $|1-z_n|<1$ and thus $\arg(z_n)\in
(-\pi/2,\pi/2)$. A direct calculation shows that then
$|\arg(z_{n+1})|<|\arg(z_n)|$ and thus, if $m>n$, then
$\arg(z_m)\in(-\pi/2,\pi/2)$. Thus by Proposition~\ref{Cor1},
(\roman{pi} and \roman{piii}), eventually $z_n\in\mathcal{V}$. We know
that $\mathcal{V}$ is a domain of analyticity of $R$. By (\ref{eq:26}),
if $H$ is analytic at $z_{n+1}=z_{n}-z_{n}^2$ then $H$ is analytic at
$z_{n}$ and by induction $H$ is analytic at $z_0$. Since $L_f$ is simply
connected, we have that $H$, and thus $R$, is analytic in $L_f$, as in
the proof of Theorem~\ref{Barrier}.

On the other hand, if we assume that $v\in\partial L$ is a periodic
point of $f$, say of period $N$, and that $R$, thus $H$, is analytic
there, relation (\ref{eq:26}) implies that $H$ is analytic at any point
on the orbit of $v$ and furthermore $H(v)=H(v)-N$, a contradiction.
Since the closure of the periodic points is $\partial L$, $\partial L$
is a singularity barrier of $H$. Furthermore $\partial L$ is in the
exterior of $\mathcal{V}$ and since $\partial \mathcal{V}$ touches the
origin tangentially along $\RR^-$, so does $\partial L$ since
$0\in\partial L$. $\Box$

\end{proof}

                         %% You may leave this command even
                                        %% if you don't have appendices.

%%      If you have an appendix, uncomment the following line and, 
%%      if desired, add a title. Then insert the text.

%%      Type body of single appendix here.

%%      ---------------------------------------------------------------------
%%      ---------------------------ACKNOWLEDGMENTS (OPTIONAL) ---------------
%%      ---------------------------------------------------------------------

%% ***** UNCOMMENT THE FOLLOWING LINE TO ADD ACKNOWLEDGMENTS.
\ack The authors are very grateful to R D Costin for many useful
discussions and comments. The authors would also like to thank R Conte,
F Fauvet, N Joshi and D Sauzin for interesting discussions.

%%      Type acknowledgments here.

%%      ---------------------------------------------------------------------
%%      --------------------------- BIBLIOGRAPHY ----------------------------
%%      ---------------------------------------------------------------------

\frenchspacing
\bibliographystyle{plain}

%%      For each reference, provide the following information:

% \bibitem{* Label *}                   %% Give a reference label.
% * Name(s) of Author(s) *              %% Enter author(s) names.
% EXAMPLE:  Gray, M., Black, F., and White, A.          

%%      Use the following template for a journal article:
% * Title of article *.                 %% Example: Existence and uniqueness
% \textit{* Abbreviated journal title *} % Example: Comm. Pure Appl. Math.
% * Volume number *:                    %% Example: 72
% * Page range *,                       %% Example: 675--690
% * Year of publication *.              %% Example: 1993
                                
%%      Use the following template for a book:
% \textit{* Title of book *}.           %% Example: Ancient Topology
% * Publisher *,                        %% Example: Wiley-Interscience
% * City of publisher *,                %% Example: New York
% * Year of publication *               %% Example: 1993

\end{document}